
\magnification=\magstephalf
\baselineskip=13pt
\hsize = 5.5truein
\hoffset = 0.5truein
\vsize = 8.5truein
\voffset = 0.2truein
\emergencystretch = 0.05\hsize

\newif\ifblackboardbold

\blackboardboldtrue


\font\titlefont=cmbx12 scaled\magstephalf
\font\sectionfont=cmbx12

\font\scriptit=cmti10 at 7pt
\font\scriptsl=cmsl10 at 7pt
\scriptfont\itfam=\scriptit
\scriptfont\slfam=\scriptsl


\newfam\bboldfam
\ifblackboardbold
\font\tenbbold=msbm10
\font\sevenbbold=msbm7
\font\fivebbold=msbm5
\textfont\bboldfam=\tenbbold
\scriptfont\bboldfam=\sevenbbold
\scriptscriptfont\bboldfam=\fivebbold
\def\bbold{\fam\bboldfam\tenbbold}
\else
\def\bbold{\bf}
\fi


\font\Arm=cmr8
\font\Ai=cmmi8
\font\Asy=cmsy8
\font\Abf=cmbx8
\font\Brm=cmr6
\font\Bi=cmmi6
\font\Bsy=cmsy6
\font\Bbf=cmbx6
\font\Crm=cmr5
\font\Ci=cmmi5
\font\Csy=cmsy5
\font\Cbf=cmbx5

\ifblackboardbold
\font\Abbold=msbm10 at 8pt
\font\Bbbold=msbm7 at 6pt
\font\Cbbold=msbm5
\fi

\def\smallmath{%
\textfont0=\Arm \scriptfont0=\Brm \scriptscriptfont0=\Crm
\textfont1=\Ai \scriptfont1=\Bi \scriptscriptfont1=\Ci
\textfont2=\Asy \scriptfont2=\Bsy \scriptscriptfont2=\Csy
\textfont\bffam=\Abf \scriptfont\bffam=\Bbf \scriptscriptfont\bffam=\Cbf
\def\rm{\fam0\Arm}\def\mit{\fam1}\def\oldstyle{\fam1\Ai}%
\def\bf{\fam\bffam\Abf}%
\ifblackboardbold
\textfont\bboldfam=\Abbold
\scriptfont\bboldfam=\Bbbold
\scriptscriptfont\bboldfam=\Cbbold
\def\bbold{\fam\bboldfam\Abbold}%
\fi
}








\newlinechar=`@
\def\forwardmsg#1#2#3{\immediate\write16{@*!*!*!* forward reference should
be: @\noexpand\forward{#1}{#2}{#3}@}}
\def\nodefmsg#1{\immediate\write16{@*!*!*!* #1 is an undefined reference@}}

\def\forwardsub#1#2{\def\newref{{#2}{#1}}}

\def\forward#1#2#3{%
\expandafter\expandafter\expandafter\forwardsub\expandafter{#3}{#2}
\expandafter\ifx\csname#1\endcsname\relax\else%
\expandafter\ifx\csname#1\endcsname\newref\else%
\forwardmsg{#1}{#2}{#3}\fi\fi%
\expandafter\let\csname#1\endcsname\newref}

\def\firstarg#1{\expandafter\argone #1}\def\argone#1#2{#1}
\def\secondarg#1{\expandafter\argtwo #1}\def\argtwo#1#2{#2}

\def\ref#1{\expandafter\ifx\csname#1\endcsname\relax
  {\nodefmsg{#1}\bf`#1'}\else
  \expandafter\firstarg\csname#1\endcsname
  ~\expandafter\secondarg\csname#1\endcsname\fi}

\def\refs#1{\expandafter\ifx\csname#1\endcsname\relax
  {\nodefmsg{#1}\bf`#1'}\else
  \expandafter\firstarg\csname #1\endcsname
  s~\expandafter\secondarg\csname#1\endcsname\fi}

\def\refn#1{\expandafter\ifx\csname#1\endcsname\relax
  {\nodefmsg{#1}\bf`#1'}\else
  \expandafter\secondarg\csname #1\endcsname\fi}



\def\widow#1{\vskip 0pt plus#1\vsize\goodbreak\vskip 0pt plus-#1\vsize}



\def\marginlabel#1{}

\def\showlabelsabove{
\font\labelfont=cmss10 at 6pt
\def\marginlabel##1{\rlap{\smash{\raise 10pt\hbox{\labelfont##1}}}}
}

\newcount\seccount
\newcount\proccount
\seccount=0
\proccount=0

\def\stdskip{\vskip 9pt plus3pt minus 3pt}
\def\stdbreak{\par\removelastskip\penalty-100\stdskip}

\def\proof{\stdbreak\noindent{\bf Proof. }}

\def\qed{\vrule height 1.2ex width .9ex depth .1ex}

\def\Box{
  \ifmmode\eqno\qed
  \else\ifvmode\removelastskip\line{\hfil\qed}
  \else\unskip\quad\hskip-\hsize
    \hbox{}\hskip\hsize minus 1em\qed\par
  \fi\stdbreak\fi}

\def\references{
  \removelastskip
  \widow{.05}
  \vskip 24pt plus 6pt minus 6 pt
  \leftline{\sectionfont References}
  \nobreak\stdskip\noindent}

\def\ifempty#1#2\endB{\ifx#1\endA}
\def\makeref#1#2#3{\ifempty#1\endA\endB\else\forward{#1}{#2}{#3}\fi}

\outer\def\section#1 #2\par{
  \removelastskip
  \global\advance\seccount by 1
  \global\proccount=0\relax
		\edef\numtoks{\number\seccount}
  \makeref{#1}{Section}{\numtoks}
  \widow{.05}
  \vskip 24pt plus 6pt minus 6 pt
  \message{#2}
  \leftline{\marginlabel{#1}\sectionfont\numtoks\quad #2}
  \nobreak\stdskip}

\def\proclamation#1#2{
  \outer\expandafter\def\csname#1\endcsname##1 ##2\par{
  \stdbreak
  \advance\proccount by 1
  \edef\numtoks{\number\seccount.\number\proccount}
  \makeref{##1}{#2}{\numtoks}
  \noindent{\marginlabel{##1}\bf #2 \numtoks\enspace}
  {\sl##2\par}
  \stdbreak}}

\def\othernumbered#1#2{
  \outer\expandafter\def\csname#1\endcsname##1{
  \stdbreak
  \advance\proccount by 1
  \edef\numtoks{\number\seccount.\number\proccount}
  \makeref{##1}{#2}{\numtoks}
  \noindent{\marginlabel{##1}\bf #2 \numtoks\enspace}}}

\proclamation{definition}{Definition}
\proclamation{lemma}{Lemma}
\proclamation{proposition}{Proposition}
\proclamation{theorem}{Theorem}
\proclamation{corollary}{Corollary}
\proclamation{conjecture}{Conjecture}

\othernumbered{example}{Example}
\othernumbered{remark}{Remark}
\othernumbered{construction}{Construction}





\input epsf

\newcount\figcount
\figcount=0
\newbox\drawing
\newcount\drawbp
\newdimen\drawx
\newdimen\drawy
\newdimen\ngap
\newdimen\sgap
\newdimen\wgap
\newdimen\egap

\def\drawbox#1#2#3{\vbox{
  \setbox\drawing=\vbox{\offinterlineskip\epsfbox{#2.eps}\kern 0pt}
  \drawbp=\epsfurx
  \advance\drawbp by-\epsfllx\relax
  \multiply\drawbp by #1
  \divide\drawbp by 100
  \drawx=\drawbp truebp
  \ifdim\drawx>\hsize\drawx=\hsize\fi
  \epsfxsize=\drawx
  \setbox\drawing=\vbox{\offinterlineskip\epsfbox{#2.eps}\kern 0pt}
  \drawx=\wd\drawing
  \drawy=\ht\drawing
  \ngap=0pt \sgap=0pt \wgap=0pt \egap=0pt
  \setbox0=\vbox{\offinterlineskip
    \box\drawing \ifgridlines\drawgrid\drawx\drawy\fi #3}
  \kern\ngap\hbox{\kern\wgap\box0\kern\egap}\kern\sgap}}

\def\draw#1#2#3{
  \setbox\drawing=\drawbox{#1}{#2}{#3}
  \advance\figcount by 1
  \goodbreak
  \midinsert
  \centerline{\ifgridlines\boxgrid\drawing\fi\box\drawing}
  \smallskip
  \vbox{\offinterlineskip
    \centerline{Figure~\number\figcount}
    \smash{\marginlabel{#2}}}
  \endinsert}

\def\nextfigtoks{%
  \advance\figcount by 1%
  \edef\numtoks{\number\figcount}%
  \advance\figcount by -1}

\def\nextfig{\nextfigtoks Figure~\numtoks}

\newif\ifgridlines
\newbox\figtbox
\newbox\figgbox
\newdimen\figtx
\newdimen\figty

\newdimen\bwd
\bwd=2sp 

\def\hline#1{\vbox{\smash{\hbox to #1{\leaders\hrule height \bwd\hfil}}}}

\def\vline#1{\hbox to 0pt{%
  \hss\vbox to #1{\leaders\vrule width \bwd\vfil}\hss}}

\def\clap#1{\hbox to 0pt{\hss#1\hss}}
\def\vclap#1{\vbox to 0pt{\offinterlineskip\vss#1\vss}}

\def\hstutter#1#2{\hbox{%
  \setbox0=\hbox{#1}%
  \hbox to #2\wd0{\leaders\box0\hfil}}}

\def\vstutter#1#2{\vbox{
  \setbox0=\vbox{\offinterlineskip #1}
  \dp0=0pt
  \vbox to #2\ht0{\leaders\box0\vfil}}}

\def\crosshairs#1#2{
  \dimen1=.002\drawx
  \dimen2=.002\drawy
  \ifdim\dimen1<\dimen2\dimen3\dimen1\else\dimen3\dimen2\fi
  \setbox1=\vclap{\vline{2\dimen3}}
  \setbox2=\clap{\hline{2\dimen3}}
  \setbox3=\hstutter{\kern\dimen1\box1}{4}
  \setbox4=\vstutter{\kern\dimen2\box2}{4}
  \setbox1=\vclap{\vline{4\dimen3}}
  \setbox2=\clap{\hline{4\dimen3}}
  \setbox5=\clap{\copy1\hstutter{\box3\kern\dimen1\box1}{6}}
  \setbox6=\vclap{\copy2\vstutter{\box4\kern\dimen2\box2}{6}}
  \setbox1=\vbox{\offinterlineskip\box5\box6}
  \smash{\vbox to #2{\hbox to #1{\hss\box1}\vss}}}

\def\boxgrid#1{\rlap{\vbox{\offinterlineskip
  \setbox0=\hline{\wd#1}
  \setbox1=\vline{\ht#1}
  \smash{\vbox to \ht#1{\offinterlineskip\copy0\vfil\box0}}
  \smash{\vbox{\hbox to \wd#1{\copy1\hfil\box1}}}}}}

\def\drawgrid#1#2{\vbox{\offinterlineskip
  \dimen0=\drawx
  \dimen1=\drawy
  \divide\dimen0 by 10
  \divide\dimen1 by 10
  \setbox0=\hline\drawx
  \setbox1=\vline\drawy
  \smash{\vbox{\offinterlineskip
    \copy0\vstutter{\kern\dimen1\box0}{10}}}
  \smash{\hbox{\copy1\hstutter{\kern\dimen0\box1}{10}}}}}

\def\figtext#1#2#3#4#5{
  \setbox\figtbox=\hbox{#5}
  \dp\figtbox=0pt
  \figtx=-#3\wd\figtbox \figty=-#4\ht\figtbox
  \advance\figtx by #1\drawx \advance\figty by #2\drawy
  \dimen0=\figtx \advance\dimen0 by\wd\figtbox \advance\dimen0 by-\drawx
  \ifdim\dimen0>\egap\global\egap=\dimen0\fi
  \dimen0=\figty \advance\dimen0 by\ht\figtbox \advance\dimen0 by-\drawy
  \ifdim\dimen0>\ngap\global\ngap=\dimen0\fi
  \dimen0=-\figtx
  \ifdim\dimen0>\wgap\global\wgap=\dimen0\fi
  \dimen0=-\figty
  \ifdim\dimen0>\sgap\global\sgap=\dimen0\fi
  \smash{\rlap{\vbox{\offinterlineskip
    \hbox{\hbox to \figtx{}\ifgridlines\boxgrid\figtbox\fi\box\figtbox}
    \vbox to \figty{}
    \ifgridlines\crosshairs{#1\drawx}{#2\drawy}\fi
    \kern 0pt}}}}

\def\nwtext#1#2#3{\figtext{#1}{#2}01{#3}}
\def\netext#1#2#3{\figtext{#1}{#2}11{#3}}
\def\swtext#1#2#3{\figtext{#1}{#2}00{#3}}
\def\setext#1#2#3{\figtext{#1}{#2}10{#3}}

\def\wtext#1#2#3{\figtext{#1}{#2}0{.5}{#3}}
\def\etext#1#2#3{\figtext{#1}{#2}1{.5}{#3}}
\def\ntext#1#2#3{\figtext{#1}{#2}{.5}1{#3}}
\def\stext#1#2#3{\figtext{#1}{#2}{.5}0{#3}}
\def\ctext#1#2#3{\figtext{#1}{#2}{.5}{.5}{#3}}


\def\hpad#1#2#3{\hbox{\kern #1\hbox{#3}\kern #2}}
\def\vpad#1#2#3{\setbox0=\hbox{#3}\dp0=0pt\vbox{\kern #1\box0\kern #2}}

\def\wpad#1#2{\hpad{#1}{0pt}{#2}}
\def\epad#1#2{\hpad{0pt}{#1}{#2}}
\def\npad#1#2{\vpad{#1}{0pt}{#2}}
\def\spad#1#2{\vpad{0pt}{#1}{#2}}

\def\nwpad#1#2#3{\npad{#1}{\wpad{#2}{#3}}}
\def\nepad#1#2#3{\npad{#1}{\epad{#2}{#3}}}
\def\swpad#1#2#3{\spad{#1}{\wpad{#2}{#3}}}
\def\sepad#1#2#3{\spad{#1}{\epad{#2}{#3}}}


\def\mimic#1#2{\setbox1=\hbox{#1}\setbox2=\hbox{#2}\ht2=\ht1\box2}


\def\stack#1#2#3{\vbox{\offinterlineskip
  \setbox2=\hbox{#2}
  \setbox3=\hbox{#3}
  \dimen0=\ifdim\wd2>\wd3\wd2\else\wd3\fi
  \hbox to \dimen0{\hss\box2\hss}
  \kern #1
  \hbox to \dimen0{\hss\box3\hss}}}


\def\sexp#1{\rlap{${}^{#1}$}}
\def\hexp#1{%
  \setbox0=\hbox{${}^{#1}$}%
  \hbox to .5\wd0{\box0\hss}}
\def\ssub#1{\rlap{${}_{#1}$}}
\def\hsub#1{%
  \setbox0=\hbox{${}_{#1}$}%
  \hbox to .5\wd0{\box0\hss}}



\def\bmatrix#1#2{{\left[\vcenter{\halign
  {&\kern#1\hfil$##\mathstrut$\kern#1\cr#2}}\right]}}

\def\rightarrowmat#1#2#3{
  \setbox1=\hbox{\smallmath\kern#2$\bmatrix{#1}{#3}$\kern#2}
  \,\vbox{\offinterlineskip\hbox to\wd1{\hfil\copy1\hfil}
    \kern 3pt\hbox to\wd1{\rightarrowfill}}\,}

\def\leftarrowmat#1#2#3{
  \setbox1=\hbox{\smallmath\kern#2$\bmatrix{#1}{#3}$\kern#2}
  \,\vbox{\offinterlineskip\hbox to\wd1{\hfil\copy1\hfil}
    \kern 3pt\hbox to\wd1{\leftarrowfill}}\,}

\def\rightarrowbox#1#2{
  \setbox1=\hbox{\kern#1\hbox{\smallmath #2}\kern#1}
  \,\vbox{\offinterlineskip\hbox to\wd1{\hfil\copy1\hfil}
    \kern 3pt\hbox to\wd1{\rightarrowfill}}\,}

\def\leftarrowbox#1#2{
  \setbox1=\hbox{\kern#1\hbox{\smallmath #2}\kern#1}
  \,\vbox{\offinterlineskip\hbox to\wd1{\hfil\copy1\hfil}
    \kern 3pt\hbox to\wd1{\leftarrowfill}}\,}






\def\bookletdims{
  \hsize=5.25truein
  \vsize=7truein
}

\def\legalbooklet#1{
  \input quire
  \bookletdims
  \htotal=7.0truein
  \vtotal=8.5truein
  \hoffset=\htotal
  \advance\hoffset by -\hsize
  \divide\hoffset by 2
  \voffset=\vtotal
  \advance\voffset by -\vsize
  \divide\voffset by 2
  \advance\voffset by -.0625truein
  \shhtotal=2\htotal
  \horigin=0.0truein
  \vorigin=0.0truein
  \shstaplewidth=0.01pt
  \shstaplelength=0.66truein
  \shthickness=0pt
  \shoutline=0pt
  \shcrop=0pt
  \shvoffset=-1.0truein
  \ifnum#1>0\quire{#1}\else\qtwopages\fi
}

\def\preview{
  \input quire
  \bookletdims
  \hoffset=0.1truein
  \vtotal=8.5truein
  \shhtotal=14truein
  \voffset=\vtotal
  \advance\voffset by -\vsize
  \divide\voffset by 2
  \advance\voffset by -.0625truein
  \htotal=2\hoffset
  \advance\htotal by \hsize
  \horigin=0.0truein
  \vorigin=0.0truein
  \shstaplewidth=0.5pt
  \shstaplelength=0.5\vtotal
  \shthickness=0pt
  \shoutline=0pt
  \shcrop=0pt
  \shvoffset=-1.0truein
  \qtwopages
}

\def\twoup{
  \input quire
  \hsize=4.79452truein 
  \vsize=7truein
  \vtotal=8.5truein
  \shhtotal=11truein
  \hoffset=-2\hsize
  \advance\hoffset by \shhtotal
  \divide\hoffset by 6
  \voffset=\vtotal
  \advance\voffset by -\vsize
  \divide\voffset by 2
  \advance\voffset by -12truept
  \htotal=2\hoffset
  \advance\htotal by \hsize
  \horigin=0.0truein
  \vorigin=0.0truein
  \shstaplewidth=0.01pt
  \shstaplelength=0pt
  \shthickness=0pt
  \shoutline=0pt
  \shcrop=0pt
  \shvoffset=-1.0truein
  \qtwopages
}


\def\enma#1{{\ifmmode#1\else$#1$\fi}}

\def\mathbb#1{{\bbold #1}}
\def\mathbf#1{{\bf #1}}

\def\NN{\enma{\mathbb{N}}}
\def\QQ{\enma{\mathbb{Q}}}
\def\RR{\enma{\mathbb{R}}}
\def\ZZ{\enma{\mathbb{Z}}}
\def\bbPP{\enma{\mathbb{P}}}


\def\aa{\enma{\mathbf{a}}}
\def\bb{\enma{\mathbf{b}}}
\def\cc{\enma{\mathbf{c}}}
\def\ee{\enma{\mathbf{e}}}

\def\uu{\enma{\mathbf{u}}}
\def\vv{\enma{\mathbf{v}}}
\def\ww{\enma{\mathbf{w}}}
\def\xx{\enma{\mathbf{x}}}

\def\zz{\enma{\mathbf{z}}}

\def\FF{\enma{\mathbf{F}}}

\def\zero{\enma{\mathbf{0}}}


\def\set#1{\enma{\{#1\}}}

\def\setthree#1#2#3{\set{#1,\, \allowbreak #2,\, \allowbreak #3}}

\def\setdef#1#2{\enma{\{\;#1\;\,|\allowbreak
  \;\,#2\;\}}}
\def\bigsetdef#1#2{\enma{\left\{\;#1\;\;\left|
  \;\;#2\;\right.\right\}}}

\def\ideal#1{\langle\,#1\,\rangle}

\def\idealfour#1#2#3#4{\ideal{#1,\, \allowbreak #2,\, \allowbreak #3,\,
  \allowbreak #4}}

\def\bigidealdef#1#2{\enma{\left\langle\;#1\;\;
  \left|\;\;#2\;\right.\right\rangle}}

\def\mtext#1{\;\,\allowbreak\hbox{#1}\allowbreak\;\,}

\def\abs#1{\enma{\left| #1 \right|}}
\def\dsum{\mathop{\hbox{$\bigoplus$}}}

\def\Tor{\mathop{\rm Tor}\nolimits}
\def\Hom{\mathop{\rm Hom}\nolimits}

\def\det{\mathop{\rm det}\nolimits}
\def\lcm{\mathop{\rm lcm}\nolimits}
\def\gcd{\mathop{\rm gcd}\nolimits}
\def\ker{\mathop{\rm ker}\nolimits}

\def\image{\mathop{\rm image}\nolimits}
\def\conv{\mathop{\rm conv}\nolimits}
\def\hull{\mathop{\rm hull}\nolimits}
\def\fiber{\mathop{\rm fiber}\nolimits}
\def\deg{\mathop{\rm deg}\nolimits}

\def\charac{\mathop{\rm char }\nolimits}
\def\sign{\mathop{\rm sign}\nolimits}
\def\Rep{\mathop{\rm Rep}\nolimits}


\forward{propD}{Proposition}{1.1}
\forward{exD}{Example}{1.9}
\forward{corK}{Corollary}{1.12}
\forward{exR}{Example}{1.10}
\forward{exA}{Example}{1.6}
\forward{rpfig}{Figure}{2}
\forward{nonfinex}{Example}{2.10}


\hbox{}
\bigskip
\centerline{\titlefont Cellular Resolutions of Monomial Modules} \bigskip
\centerline{\bf Dave Bayer \quad Bernd Sturmfels}
\bigskip


\vskip .3cm

\noindent {\sl {\bf Abstract:}
We construct a canonical free resolution
for arbitrary  monomial modules and lattice ideals.
This includes monomial ideals and  defining ideals of toric varieties,
and it generalizes our joint results with Irena Peeva
for generic ideals [BPS],[PS].}

\vskip .6cm

\centerline{\bf Introduction}

\vskip .2cm

\noindent
Given a field $k$, we consider the Laurent polynomial
ring $T = k[x_1^{\pm 1}\!\!,\, \dots ,\, x_n^{\pm 1}]$
as a module over the polynomial ring $S = k[x_1,\,\dots,\, x_n]$.
The module structure comes from the natural inclusion of semigroup algebras
$S = k[\NN^n] \, \subset \,  k[\ZZ^n]= T$.
A {\it monomial module} is an $S$-submodule of $T$ which is
generated by  monomials $\, \xx^\aa = x_1^{a_1} \cdots x_n^{a_n}$,
$\aa \in \ZZ^n$. Every monomial module $M$ has a unique
minimal generating set of monomials.
Of special interest are the two cases when this generating set
is finite or when it forms a group under multiplication.
In the first case $M$ is isomorphic to a {\it monomial ideal} in $S$.
In the second case $M$ coincides with the {\it lattice module}
 $$ M_L \quad := \quad
S \, \{ \xx^\aa \,\,|\,\, \aa \in L \} \quad = \quad
k \, \{ \xx^\bb \,\,|\,\, \bb \in \NN^n+L \} \quad \subset \quad T.$$
for some sublattice $L\subset\ZZ^n$ whose intersection
with $\NN^n$ is the origin $\zero=(0,\ldots,0)$.

We shall derive free resolutions of $M$
from regular cell complexes whose vertices
are the generators of $M$
and whose faces are labeled by the
least common multiples of their vertices.
The basic theory of such cellular resolutions
is developed in Section~1.

Our main result is the construction of the {\it hull resolution}
in Section~2.  We rescale the
exponents of the monomials in $M$, so that their convex hull in
$\RR^n$ is a polyhedron $P_t$ whose bounded faces support a free
resolution of $M$. This resolution is new and interesting
even for monomial ideals. It need not be
minimal, but, unlike minimal resolutions, it
respects symmetry and is free from arbitrary choices.

In Section~3 we relate the lattice module $M_L$ to the
 $\ZZ^n/L$-graded {\it lattice ideal\/}
$$I_L \quad = \quad \bigidealdef{\xx^\aa-\xx^\bb}{\aa-\bb\in L} \quad
\subset \quad S.$$
This class of ideals includes ideals defining {\it toric varieties}.
We express the cyclic $S$-module $S/I_L$ as the quotient of the
infinitely generated $S$-module  $M_L$ by the action of $L$.
In fact, we like to think of $M_L$ as the ``universal cover'' of $I_L$.
Many questions about $I_L$ can thus be reduced to
questions about $M_L$. In particular, we obtain the hull resolution
of a lattice ideal $I_L$ by taking the hull resolution of $M_L$ modulo $L$.

This paper is inspired by the work of
Barany, Howe and Scarf [BHS] who introduced
the polyhedron $P_t$ in the context of  integer programming.
The hull resolution generalizes results in [BPS] for
generic monomial ideals and in [PS] for generic lattice ideals.
In these generic cases the hull resolution is minimal.
We remark that deriving algebraic syzygies from polytopes is
a natural operation also in other contexts such as
associative algebras [AB] and algebraic K-theory [KS].

\section{} Cellular resolutions

\noindent Fix a subset $\setdef{\aa_j}{j\in I} \subset
\ZZ^n\!$, for an index set $I$ which need not be finite,
and let $M  $ be the monomial module
generated by the monomials $m_j = \xx^{\aa_j}$, $ j \in I$.
Let $X$ be a {\it regular cell complex\/} having $I$ as its set
of vertices, and equipped with a choice of an
{\it incidence function} $\varepsilon(F,F')$ on pairs of faces.
We recall from [BH, Section 6.2] that $\,\varepsilon\,$
takes values in \set{0,1,-1}, that
$\varepsilon(F,F') = 0$ unless $F'$ is a facet of $F$, that
$\varepsilon(\{j\},\emptyset) =1$ for all vertices $j\in I$, and
that for any codimension 2 face $F'$ of $F$,
$$\varepsilon(F,F_1) \varepsilon(F_1,F') \,\, + \,\,
\varepsilon(F,F_2) \varepsilon(F_2,F') \quad = \quad 0$$
where $F_1$, $F_2$ are the two facets of $F$ containing $F'$.
The prototype of a regular cell complex is the set of faces of a
convex polytope. The incidence function $\varepsilon$ defines a differential
$\partial$ which can be used to compute the homology of $X$. Define
the {\it augmented oriented chain complex} $\widetilde{C}(X;k) =
\dsum_{F\in X}\; kF$, with differential
 $$ \partial F \quad = \quad \sum_{F'\in X} \, \varepsilon(F,F') \, F'.$$
 The {\it reduced cellular homology group} $\widetilde{H}_i(X;k)$ is
the $i$th homology of $\widetilde{C}(X;k)$, where faces of $X$
are indexed by their dimension. The {\it oriented
chain complex} $C(X;k) = \dsum_{F\in X,\, F\ne\emptyset}\; kF$ is
obtained from  $\widetilde{C}(X;k)$  by dropping the
contribution of the empty face.
It computes the ordinary homology groups $H_i(X;k)$ of $X$.

The cell complex $X$ inherits a $\ZZ^n$-grading from the
generators of $M$ as follows. Let $F$ be a nonempty face of $X$.
We identify $F$ with its set of vertices, a finite subset of $I$.
Set $m_F := \lcm \setdef{m_j}{j\in F}$. The exponent vector of
the monomial $m_F$ is the {\it join}
$ \, \aa_F := \bigvee\setdef{\aa_j}{j\in F}$
in $\ZZ^n $. We call  $\aa_F $
the {\it degree} of the face $F$.

Homogenizing the differential $\partial$ of $C(X;k)$
yields a $\ZZ^n$-graded chain complex of $S$-modules. Let $SF$ be the free
$S$-module with one generator $F$ in degree $\aa_F$. The {\it
cellular complex} $\, \FF_X \,$ is the $\ZZ^n$-graded
$S$-module $\dsum_{F\in X,\, F\ne\emptyset}  SF$ with differential
 $$ \partial \, F \quad = \quad \sum_{F'\in X,\, F'\ne\emptyset} \;
\varepsilon(F,F') \; {m_F \over m_{F'}} \; F'.$$
 The homological degree of each face $F$ of $X$ is its dimension.

For each degree $\bb\in\ZZ^n$, let $X_{\preceq\bb}$ be the
subcomplex of $X$ on the vertices of degree $\preceq\bb$, and let
$X_{\prec\bb}$ be the subcomplex of $X_{\preceq\bb}$ obtained by
deleting the faces of degree \bb.
For example, if there is a unique vertex $j$ of degree $\preceq\bb$, and
$\aa_j=\bb$, then $X_{\preceq\bb} = \set{\set{j},\emptyset}$ and
$X_{\prec\bb}=\set{\emptyset}$. A full subcomplex on no vertices is the
acyclic complex \set{}, so if there are no vertices of degree
$\preceq\bb$, then $X_{\preceq\bb} = X_{\prec\bb} = \set{}$.

 The following proposition generalizes [BPS,~Lemma~2.1] to cell
complexes:

\proposition{propD} 
 The complex $\FF_X$ is a free resolution of $M$ if and only if
$X_{\preceq\bb}$ is acyclic over $k$ for all degrees $\bb$.
In this case we call $\FF_X$ a {\it cellular resolution} of $M$.

\proof
The complex $\FF_X$ is $\ZZ^n$-graded. The degree \bb\ part of $\FF_X$ is
precisely the oriented chain complex $C(X_{\preceq\bb};k)$. Hence $\FF_X$ is a
free resolution of $M$ if and only if $H_0(X_{\preceq\bb};k) \cong k$ for
$\xx^\bb \in M$, and otherwise $H_i(X_{\preceq\bb};k) = 0$ for all $i$ and
all $\bb$. This condition is equivalent to
$\widetilde{H}_i(X_{\preceq\bb};k) = 0$  for all $i$ (since $\xx^\bb \in M
\,$ if and only if $ \emptyset \in X_{\preceq\bb}$) and thus  to
$X_{\preceq\bb}$ being acyclic. \Box

\remark{remU} 
Fix $\bb\in\ZZ^n$. The set of generators of $M$ of degree
$\preceq\bb$ is finite. It generates a monomial module $M_{\preceq \bb}$
isomorphic  to an ideal (up to a degree shift).

\corollary{makefinite}
The cellular complex $\FF_X$ is a resolution of  $M$ if and only if the
cellular complex $\FF_{X_{\preceq\bb}}$ is a resolution of the monomial
ideal $M_{\preceq\bb}$ for all $\bb \in \ZZ^n$.

\proof
 This follows from \ref{propD} and the identity
$\,(X_{\preceq\bb})_{\preceq\cc} = X_{\preceq \, \bb \wedge\cc}$.
 \Box

\remark{distinct}
A cellular resolution $\FF_X$ is a minimal resolution if and only
if any two comparable faces $F' \subset F$ of the complex $X$
have distinct degrees $\,\aa_F \not= \aa_{F'}$.

\vskip .2cm

The simplest example of a cellular resolution is
the {\it Taylor resolution\/} for monomial ideals  [Tay].
The  Taylor resolution is easily generalized to
arbitrary monomial modules $M$ as follows.
Let $\setdef{m_j}{j \in I}$ be the
minimal generating set of $M$. Define $\Delta$ to be
the simplicial complex  consisting of all finite subsets
of $I$, equipped with the standard incidence function $\varepsilon(F,F') =
(-1)^j$ if $F\setminus F'$ consists of the $j$th element of $F$. The Taylor
complex of $M$ is the cellular complex $\FF_\Delta$.

\proposition{propT} 
 The Taylor complex $\FF_\Delta$ is a resolution of $M$.

\proof
 By \ref{propD} we need to show that each subcomplex $\Delta_{\preceq\bb}$
of $\Delta$ is acyclic. $\Delta_{\preceq\bb}$ is the full simplex on the
set of vertices \setdef{j\in I}{\aa_j \preceq \bb}.
This set is finite by \ref{remU}.
Hence $\Delta_{\preceq\bb}$ is a finite simplex, which is acyclic.
 \Box

The Taylor resolution $\FF_\Delta$
 is typically far from minimal. If $M$ is
infinitely generated, then $\Delta$ has faces of arbitrary dimension
and $\FF_\Delta$  has infinite length.
Following [BPS, \S 2] we note that
every simplicial complex $X\subset\Delta$
defines a  submodule $ \FF_X \subset \FF_\Delta$ which is closed under the
differential $\partial$. We call $\FF_X$ the {\it restricted Taylor
complex\/} supported on $X$. $\FF_X$ is a resolution of $M$ if and only
if the hypothesis of \ref{propD} holds, with cellular homology
specializing to simplicial homology.

\example{exA} 
Consider the monomial ideal
$M=\idealfour{a^2b}{ac}{b^2\!}{bc^2}$ in $S=k[a,b,c]$.
\nextfig\ shows a truncated ``staircase diagram'' of
unit cubes representing the monomials in $S \backslash M$, and shows
two simplicial complexes $X$ and $Y$ on the generators of $M$.
Both are two triangles sharing an edge.
Each vertex, edge or triangle is labeled by its degree.
The notation {\bf 210}, for example,
represents the degree $(2,1,0)$ of $a^2b$.

\draw{70}{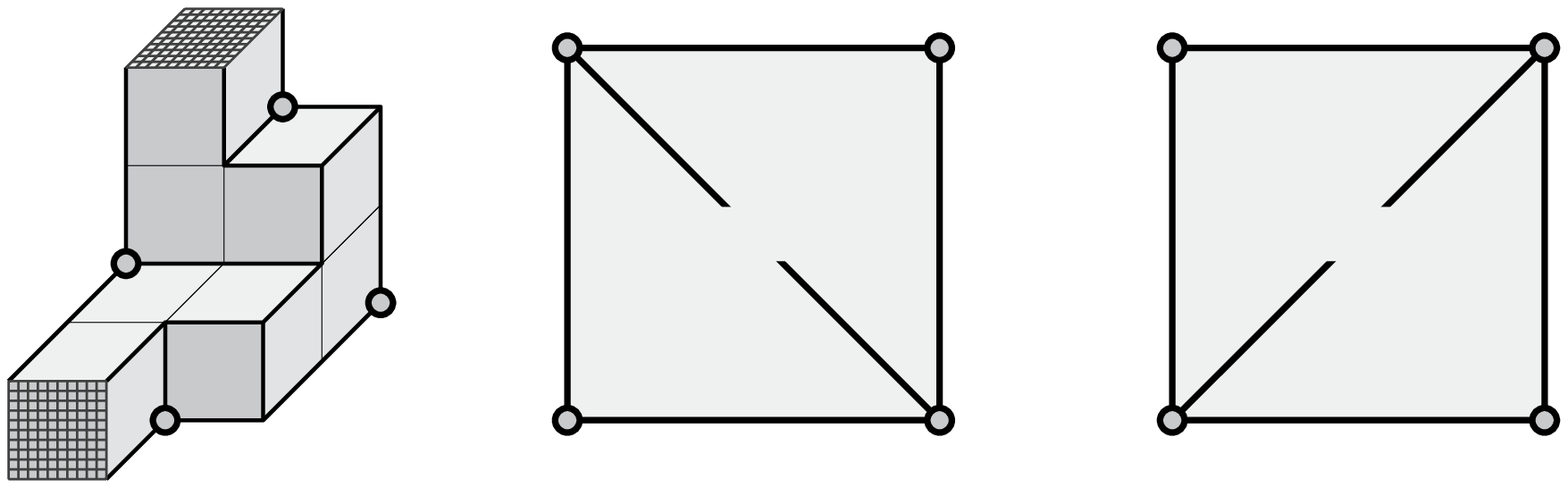}{
 \setext{.0685}{.475}{\sepad{3pt}{3pt}{$ac$}}
 \nwtext{.107}{.1}{\nwpad{1pt}{1pt}{$a^2b$}}
 \nwtext{.248}{.352}{\wpad{2pt}{$b^2$}}
 \swtext{.185}{.81}{\swpad{3pt}{3pt}{$bc^2$}}
 \setext{.353}{.94}{\sepad{2pt}{2pt}{\stack{4pt}{$ac$}{\sevenbf 101}}}
 \etext{.358}{.518}{\epad{3pt}{\sevenbf 211}}
 \netext{.353}{.105}{\nepad{2pt}{2pt}{\stack{4pt}{\sevenbf 210}{$a^2b$}}}
 \stext{.48}{.919}{\spad{3pt}{\sevenbf 112}}
 \ctext{.48}{.518}{\sevenbf 121}
 \ntext{.48}{.122}{\npad{3pt}{\sevenbf 220}}
 \swtext{.605}{.94}{\swpad{2pt}{2pt}{\stack{4pt}{$bc\hexp 2$}{\sevenbf
   012}}}
 \wtext{.602}{.518}{\wpad{3pt}{\sevenbf 022}}
 \nwtext{.605}{.105}{\nwpad{2pt}{2pt}{\stack{4pt}{\sevenbf 020}{$b\hexp
   2$}}}
 \ctext{.438}{.38}{\sevenbf 221}
 \ctext{.522}{.659}{\sevenbf 122}
 \setext{.742}{.94}{\sepad{2pt}{2pt}{\stack{4pt}{$ac$}{\sevenbf 101}}}
 \etext{.747}{.518}{\epad{3pt}{\sevenbf 211}}
 \netext{.742}{.105}{\nepad{2pt}{2pt}{\stack{4pt}{\sevenbf 210}{$a^2b$}}}
 \stext{.869}{.919}{\spad{3pt}{\sevenbf 112}}
 \ctext{.869}{.518}{\sevenbf 212}
 \ntext{.869}{.122}{\npad{3pt}{\sevenbf 220}}
 \swtext{.994}{.94}{\swpad{2pt}{2pt}{\stack{4pt}{$bc\hexp 2$}{\sevenbf
   012}}}
 \wtext{.991}{.518}{\wpad{3pt}{\sevenbf 022}}
 \nwtext{.994}{.105}{\nwpad{2pt}{2pt}{\stack{4pt}{\sevenbf 020}{$b\hexp
   2$}}}
 \ctext{.827}{.659}{\sevenbf 212}
 \ctext{.911}{.38}{\sevenbf 222}
 \ntext{.13}{0}{\npad{10pt}{$M$}}
 \ntext{.48}{0}{\npad{10pt}{$X$}}
 \ntext{.869}{0}{\npad{10pt}{$Y$}}
 }

\noindent
By \ref{propD},
the complex $X$ supports the minimal free resolution $\FF_X =$
 $$ 0 \rightarrow
 S^2 \rightarrowmat{2pt}{4pt}{\!-b&0\cr c&0\cr 0&-b\cr \!-a&c\cr
     0&-a&\cr}
 S^5 \rightarrowmat{2pt}{4pt}{c&b&0&0&0\cr
     \!-ab&0&bc&b\sexp2&0\cr
     0&0& \,-a&0&b\cr
     0&-a\sexp2&0&-ac&\,-c\sexp2&\,\cr}
 S^4 \rightarrowmat{3pt}{6pt}{a^2b & ac  & bc^2 & b^2  \cr}
 M \rightarrow 0.
 $$
The complex $Y$ fails the criterion of \ref{propD},
and hence $\FF_Y$ is not exact: if $\bb=(1,2,1)$ then
$Y_{\preceq\bb}$ consists of the two vertices $ac$ and $b^2$, and is not
acyclic. \Box

We next present four examples which are not restricted Taylor complexes.

\example{exBH}  
Let $M$ be a {\it Gorenstein ideal of height $3$}
generated by $m$ monomials.
It is shown in  [BH1, \S 6] that
 the minimal free resolution
of $M$ is the cellular resolution $\, \FF_X \,:\,
 0 \rightarrow S \rightarrow S^m \rightarrow S^m \rightarrow S
\rightarrow 0 \,$ supported on a {\it convex $m$-gon}.

\example{excogen}  %
A monomial ideal  $M$ is  {\it co-generic} if its
no variable occurs to the same power in two distinct
irreducible components
$\,\langle x_{i_1}^{r_1},
 x_{i_2}^{r_2}, \ldots, x_{i_s}^{r_s} \rangle \,$ of $M$. It is shown in
[Stu2] that the minimal resolution of a co-generic
monomial ideal is a cellular resolution
$\FF_X$ where $X$ is the complex of bounded faces
of a {\it simple polyhedron}.

\example{exD}
Let $u_1,\ldots,u_n$ be distinct integers and $M$ the module generated by
the $\, n\, !\, $ monomials $\,
x_{\pi(1)}^{u_1} x_{\pi(2)}^{u_2} \cdots
x_{\pi(n)}^{u_n} \,$ where  $\pi$ runs over all permutations of
$\,\{1,2,\ldots,n \} $.  Let $X$ be the complex of all faces of the
{\it permutohedron} [Zie, Example 0.10], which is the convex hull
of the $\,n \,! \,$ vectors $\bigl(\pi(1),,\ldots,\pi(n)\bigr)$ in $\RR^n$.
It is known [BLSWZ, Exercise~2.9] that the
$i$-faces $F$ of $X$ are indexed by chains
$$ \emptyset \;=\; A_0 \subset A_1 \subset \;\;\ldots\;\; \subset
A_{n-i-1} \subset A_{n-i} \;=\; \{u_1,u_2,\ldots,u_n\} . $$
We assign the following
monomial degree to the
$i$-face $F$ indexed by this chain:
$$ \xx^{\aa_F} \quad = \quad
\prod_{j=1}^{n-i} \; \prod_{r \in  A_j \backslash A_{j-1}}  \!\!
x_r^{\max\{ A_j \backslash A_{j-1}\}}. $$
It can be checked (using our results in \S 2)
that the conditions in \ref{propD} and \ref{distinct}
are satisfied.
Hence $\FF_X$  is the minimal free resolution of $M$. \Box

\example{exR}
 \def\rptwo{\enma{\RR\bbPP^2}}
Let $S=k[a,b,c,d,e,f]$. Following [BH, page 228] we consider the
Stanley-Reisner ideal of the  minimal triangulation
of the  {\it real projective plane} \rptwo,
$$M \;=\; \ideal{abc,\; abf,\; ace,\; ade,\; adf,\; bcd,\; bde,\;
bef,\; cdf,\; cef}.$$
The dual in \rptwo\ of this triangulation
is a cell complex $X$ consisting of six pentagons.
The ten vertices of $X$ are labeled by the
generators of $M$. We illustrate $ X \simeq $ \rptwo \ as the disk shown
on the left in \nextfig; antipodal points on the boundary are to be
identified. The small pictures on the right will be
discussed in Example 2.14.

\draw{70}{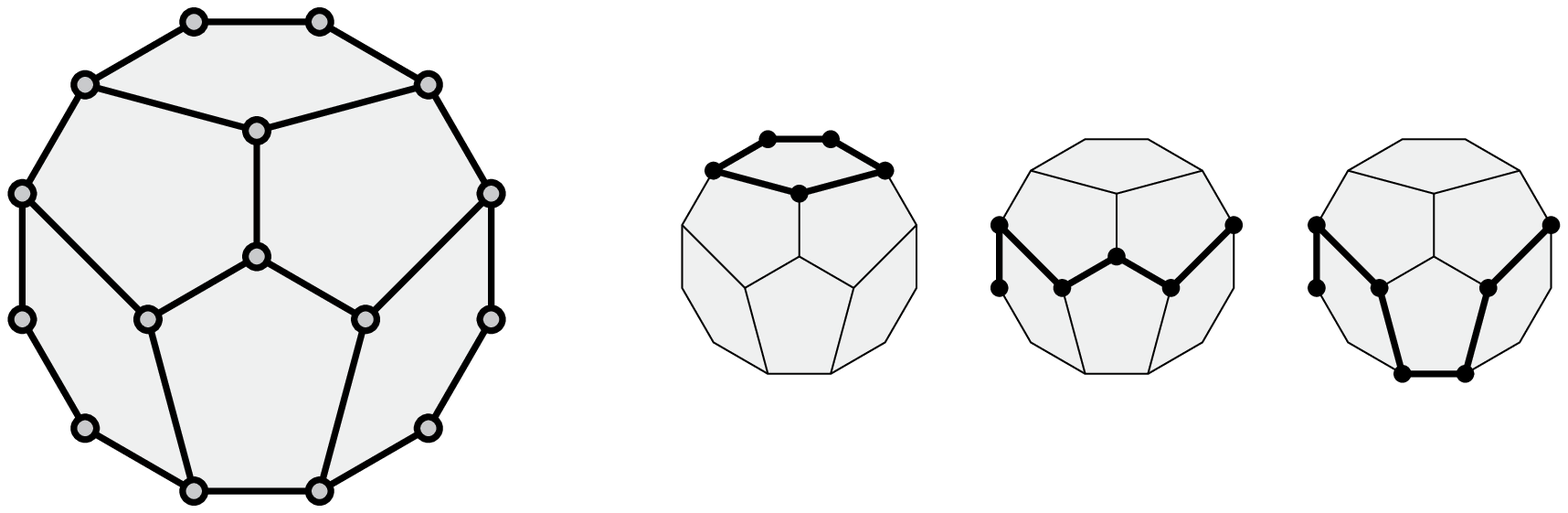}{
 \ntext{.161}{.475}{\npad{2pt}{$abc$}}
 \stext{.0915}{.406}{\swpad{4pt}{3pt}{$ace$}}
 \stext{.233}{.406}{\sepad{4pt}{3pt}{$abf$}}
 \nwtext{.169}{.728}{\npad{1pt}{$bcd$}}
 \stext{.12}{1}{\swpad{2pt}{2pt}{$bef$}}
 \stext{.202}{1}{\swpad{2pt}{2pt}{$cef$}}
 \swtext{.28}{.866}{\swpad{2pt}{1pt}{$cdf$}}
 \wtext{.322}{.626}{\wpad{2pt}{\mimic{$a$}{$adf$}}}
 \wtext{.322}{.373}{\wpad{2pt}{\mimic{$a$}{$ade$}}}
 \nwtext{.278}{.132}{\nwpad{2pt}{1pt}{$bde$}}
 \ntext{.202}{.006}{\nwpad{1pt}{2pt}{$bef$}}
 \ntext{.12}{.006}{\nwpad{1pt}{2pt}{$cef$}}
 \netext{.043}{.14}{\nepad{2pt}{1pt}{$cdf$}}
 \etext{.0015}{.373}{\epad{2pt}{\mimic{$a$}{$adf$}}}
 \etext{.0015}{.626}{\epad{2pt}{\mimic{$a$}{$ade$}}}
 \setext{.04}{.85}{\sepad{2pt}{1pt}{$bde$}}
 \ctext{.161}{.88}{$\widehat{a}$}
 \ctext{.061}{.31}{$\widehat{b}$}
 \ctext{.261}{.31}{\mimic{$\widehat{b}$}{$\widehat{c}$}}
 \ctext{.161}{.25}{$\widehat{d}$}
 \ctext{.241}{.63}{\mimic{$\widehat{f}$}{$\widehat{e}$}}
 \ctext{.081}{.63}{$\widehat{f}$}
 \ntext{.51}{.264}{\npad{10pt}{\stack{6pt}{$a\,=\,0$}{6 cycles}}}
 \ntext{.714}{.264}{\npad{10pt}{\stack{6pt}{$a\,=\,1$}{6 cycles}}}
 \ntext{.918}{.264}{\nwpad{10pt}{2pt}{\stack{6pt}{$b+c+d\,=\,1$}{10
cycles}}}
 }

\noindent
If $\charac k \ne 2$ then $X$ is acyclic over $k$ and the cellular complex
$\FF_X$ coincides with  the minimal free resolution
$\, 0 \rightarrow S^6 \rightarrow S^{15} \rightarrow S^{10} \rightarrow M $.
If  $\charac k = 2$ then $X$ is not acyclic over $k$, and the cellular
complex $\FF_X$ is not a resolution of $M$. \Box

Returning to the general theory, we next present a formula for
the {\it Betti number} $\beta_{i,\bb} = \dim\Tor_i(M,k)_\bb \,$ which
is the number of minimal $i$th syzygies  in degree $\bb$.
The degree $\bb \in \ZZ^n$ is called a {\it Betti degree} of $M$
if $\beta_{i,\bb}\not= 0$ for some $i$.

\theorem{propA} 
 If $\, \FF_X$ is a cellular resolution of a monomial module $M$ then
 $$\beta_{i,\bb} \quad = \quad \dim H_i(X_{\preceq\bb}, X_{\prec\bb};k)
\quad = \quad \dim \widetilde{H}_{i-1} (X_{\prec\bb};k),$$
 where $H_*$ denotes relative homology and $\widetilde{H}_*$ denotes
reduced homology.

\proof
We compute $\Tor_i(M,k)_\bb$ as the $i$th homology of the complex
of vector spaces $(\FF_X\otimes_S k)_\bb$. This complex equals
the chain complex
$\widetilde{C}(X_{\preceq\bb},X_{\prec\bb};k)$
which computes the relative homology with coefficients in $k$ of
the pair $(X_{\preceq\bb},X_{\prec\bb})$. Thus
 $$\Tor_i(M,k)_\bb \quad = \quad
H_i(X_{\preceq\bb},X_{\prec\bb};k).$$
Since $X_{\preceq\bb}$ is acyclic,
the long exact sequence of homology groups looks like
 $$ 0 \,=\, \widetilde{H}_i(X_{\preceq\bb}; k) \,\rightarrow
\, H_i(X_{\preceq\bb},X_{\prec\bb}; k)
\, \rightarrow \, \widetilde{H}_{i-1}(X_{\prec\bb}; k)
\, \rightarrow \, \widetilde{H}_{i-1}(X_{\preceq\bb}; k) \,=\, 0 .$$
We conclude that the two vector spaces in the middle are isomorphic.
 \Box

A subset $Q\subset\ZZ^n$ is an {\it order ideal} if $\bb\in Q$
and  $\cc \in \NN^n$ implies $\bb-\cc\in Q$. For a
$\ZZ^n$-graded cell complex $X$ and an order
ideal $Q$ we define the {\it order ideal complex\/}
$\,X_Q \, = \, \bigsetdef{F\in X}{\aa_F\in Q}$.
Note that $\,X_{\prec \bb}$ and $\,X_{\preceq \bb}$ are
special cases of this.

\corollary{corK} 
 If $\FF_X$ is a cellular resolution of $M$ and
 $Q\subset\ZZ^n$ an order ideal which contains the
 Betti degrees of $M$,
then $\FF_{X_Q}$ is also a cellular resolution of $M$.

\proof
By \ref{makefinite} and the identity
$\,(X_Q)_{\preceq \bb} =  (X_{\preceq \bb})_Q$, it suffices
to prove this for the case where $M$ is a monomial ideal and $X$ is finite.
We proceed by induction on the number of faces in $X \backslash X_Q$.
If $X_Q = X$ there is nothing to prove.  Otherwise
pick $\cc \in \ZZ^n \backslash Q$ such that $X_{\preceq \cc} = X$ and
$X_{\prec \cc} \not= X$. Since $\cc$ is not a
Betti degree, \ref{propA} implies that the complex $X_{\prec \cc}$ is acyclic.
For any $\bb \in \ZZ^n $, the complex
$\,(X_{\prec \cc})_{\preceq \bb}\,$ equals either
$X_{\prec \cc}$ or $X_{\preceq\,  \bb \wedge \cc}$
and is hence acyclic. At this point we replace $X$ by the
proper subcomplex $X_{\prec \cc}$, and we are done by induction.
\Box

By \ref{propT} and \ref{propA}, the Betti  numbers $\beta_{i,\bb}$ of $M$
are given by the reduced homology of $\Delta_{\prec\bb}$. Let us compare
that formula for $\beta_{i,\bb}$ with the following formula which is due
independently to Hochster [Ho] and Rosenknop [Ros].

\corollary{propH} 
 The Betti numbers of $M$ satisfy
 $ \,\beta_{i,\bb} \, = \,\dim \widetilde{H}_i(K_\bb;k)$
where $K_\bb$ is the simplicial complex
$\{\, \sigma \subseteq\set{1,\ldots,n}\,|\, M \! \mtext{has a generator of
degree} \! \!\preceq\bb- \sigma\}$. Here each face $\sigma$
of $K_\bb$ is identified with its
characteristic vector in $\{0,1\}^n$.

\proof
 For $i \in \{1,\ldots,n\}$ consider the subcomplex of
$\Delta_{\prec \bb}$ consisting of all faces $F$ with degree $\aa_F
\preceq \bb - \{i\}$. This subcomplex is a full simplex. Clearly, these $n$
simplices cover $\Delta_{\prec \bb}$. The nerve of this cover by
contractible subsets is the simplicial complex $K_\bb$. Therefore, $K_\bb$
has the same reduced homology as $\Delta_{\prec \bb}$.
 \Box

\section{} The hull resolution

\noindent
Let $M$ be a monomial module in $T = k[x_1^{\pm 1}\!\!,\, \dots ,\,
x_n^{\pm 1}]$. In this section we apply convexity methods to construct a
canonical cellular resolution of $M$. For $\aa \in \ZZ^n$ and $t \in \RR$
we abbreviate $t^\aa = (t^{a_1}, \dots, t^{a_n})$. Fix
any real number $t $ larger than $ (n+1) \,!\,
= 2 \cdot 3 \cdot \cdots \cdot (n+1)$.  We define
$P_t$ be the convex hull of the point set
$$\setdef{ t^\aa }{\aa \mtext{is the exponent of a monomial}
\xx^\aa\in M}
\quad \subset \quad \RR^n.$$
The set $P_t$ is a closed, unbounded  $n$-dimensional convex polyhedron.

\proposition{propO}
The vertices of the polyhedron  $P_t$ are precisely the points
$ t^\aa = (t^{a_1},\ldots,t^{a_n})  $
for which the monomial $ \xx^\aa
= x_1^{a_1}\! \cdots x_n^{a_n} $ is a minimal generator of $M$.

\proof
Suppose $\xx^\aa \in M$ is not a minimal generator of $M$.
Then $M$ contains both  $\xx^{\aa+\ee_i} = \xx^\aa x_i$
and  $\xx^{\aa-\ee_i} = \xx^\aa /x_i$ for some $i$.
The line segment $\conv \{t^{\aa-\ee_i},t^{\aa+\ee_i} \}$ lies in $P_t$
and contains $t^\aa$ in its relative interior. Therefore $t^\aa$
is not a vertex of $P_t$.

Next, suppose $\xx^\aa \in M$ is a minimal generator of $M$.
Let $\vv = t^{-\aa}$, so $\vv\cdot t^{\aa} = n$. For any other
exponent $\bb$ of a monomial in $M$, we have $b_i \ge a_{i}+1$ for some
$i$, so $$ \vv \cdot t^\bb \quad = \quad
\sum_{j=1}^n t^{b_j - a_j} \quad \geq
\quad t^{b_i-a_{i}} \quad \ge \quad t \quad > \quad (n+1) \,!
\quad > \quad n .$$
 Thus, the inner normal vector $\vv$ supports
$t^{\aa}$ as a vertex of $P_t$.
\Box

\corollary{}
$\, P_t \,\,\, = \,\,\, \RR_+^n \,\, + \,\,\conv\,
\setdef{ t^\aa }{\xx^\aa \mtext{is  a minimal generator}
\xx^\aa \mtext{of} M}$.

Our first goal is to establish the following combinatorial result.

\theorem{indepthm} The face poset of the
polyhedron $P_t$ is independent of $t$ for $\,t > (n+1) \,!$.
The same holds for the subposet of all bounded faces of $P_t$.

\proof
The face poset of $P_t$ can be computed
as follows. Let $C_t \subset \RR^{n+1}$  be the cone spanned
by the vectors $\,(t^\aa ,1 ) \,$ for all minimal generators
$\xx^\aa$ of $M$ together with the unit vectors   $\,(\ee_i, 0)\,$
for $i=1,\ldots,n$. The faces of $P_t$ are in order-preserving
bijection with the faces of $C_t$ which do not lie in the hyperplane ``at
infinity'' $\,x_{n+1} = 0$. A face of $P_t$ is bounded if and only if the
corresponding face of $C_t$ contains none of the vectors $(\ee_i,0)$.
It suffices to prove that the
face poset of $C_t$ is independent of $t$.

For any $(n+1)$-tuple of generators of $C_t$ consider
the sign of their determinant
$$
\sign \,\,
\det \bmatrix{2pt}{
\ee\ssub{i_0} & \;\;\cdots  & \ee\ssub{i_r} & \;\;\;
t^{\aa\ssub{j_1}} &
\;\;\cdots & t^{\aa\ssub{j_{n-r}}} & \;\quad\cr
   0   & \cdots &   0     &    1      & \cdots &   1 \cr}
\quad \in \quad \{-1,0,+1 \}. \eqno (1) $$
The list of
these signs forms the (possibly infinite) {\it oriented matroid}
associated with the cone $C_t$.  It is known
(see e.g.~[BLWSZ]) that the  face poset of $C_t$
is determined by its oriented matroid.
It therefore suffices to show that the sign
of the determinant in (1) is independent of $t$ for $t > (n+1) \, !$.
This follows from the next lemma.
\Box

\lemma{factoriallem}
Let $a_{ij}$ be integers for $\,1 \leq i,j \leq r$.
Then the Laurent polynomial
$ f(t) = \det \bigl( (t^{a_{ij}})_{1 \leq i,j \leq r} ) \,$
either vanishes identically or
has no real roots for $t > r \,!$.

\proof
Suppose that $f$ is not zero and write
$\,f(t) = c_\alpha t^\alpha + \sum_\beta c_\beta t^\beta$,
where the first term has  the highest degree in $t$. For $t > r!$
we have the chain of inequalities
$$ | \sum_\beta c_\beta \cdot t^\beta |
  \, \leq \,
 \sum_\beta  |c_\beta | \cdot t^\beta
  \, \leq \,
 (\sum_\beta  |c_\beta | ) \cdot t^{\alpha-1}
  \, < \,  r \, ! \cdot  t^{\alpha-1}
  \, < \, t^\alpha
  \, \le \, |c_\alpha \cdot t^\alpha| .$$
Therefore $f(t)$ is nonzero, and
$\, \sign\bigl(f(t)\bigr) = \sign(c_\alpha)$. \Box

In the proof of \ref{indepthm} we  are using \ref{factoriallem} for $r=n+1$.
Lev Borisov and Sorin Popescu constructed examples
of matrices which show
that the exponential lower bound for $t$ is necessary in \ref{factoriallem},
and also in \ref{indepthm}.

We are now ready to define the hull resolution
and state our main result.
The {\it hull complex} of a monomial module $M$,
denoted $\hull(M)$, is the complex of bounded
faces of the polyhedron $P_t$ for large $t$.
\ref{indepthm} ensures that $\hull(M)$ is well-defined
and depends only on $M$. The vertices of $\hull(M)$
are labeled by the generators of $M$, by \ref{propO},
and hence the complex $\,\hull(M) \,$ is $\ZZ^n$-graded.
Let $\FF_{\hull(M)}$ be the complex of free $S$-modules
derived from $\hull(M)$ as in Section~1.

\theorem{thmH}
 The cellular complex $\FF_{\hull(M)}$ is a free resolution of $M$.

\proof
 Let $X = (\hull(M))_{\preceq\bb}$ for some degree \bb; by \ref{propD} we
need to show that $X$ is acyclic. This is immediate if $X$ is empty or a
single vertex. Otherwise choose $\,t > (n+1)\,! \,$
 and let $\vv=t^{-\bb}$. If $t^\aa$ is a vertex of $X$ then $\aa \prec\bb$,
so
 $$\vv \cdot t^{\aa} \quad =  \quad t^{-\bb} \cdot t^{\aa} \quad < \quad
t^{-\bb} \cdot  t^\bb \quad = \quad n,$$
 while for any other $\xx^\cc \in M$ we have $c_{i}\ge b_i+1$ for some $i$,
so
 $$\vv \cdot t^{\aa} \quad = \quad t^{-\bb} \cdot t^{\cc} \quad \ge
\quad t^{c_i-b_i} \quad \ge \quad
 t \quad > \quad n.$$
Thus, the hyperplane $H$ defined by $\vv\cdot\xx = n$ separates the
vertices of $X$ from the remaining vertices of $P_t$.
Make a projective transformation which moves $H$ to infinity. This
expresses $X$ as the complex of bounded faces of a convex polyhedron, a
complex which is known to be contractible, e.g.~[BLSWZ,
Exercise 4.27 (a)].
\Box

We call $\FF_{\hull(M)}$ the {\it hull resolution} of $M$.
Let us see that the hull resolution generalizes the {\it Scarf
complex} introduced in [BPS]. This is the
simplical complex
 $$\Delta_M \quad = \quad \setdef{F \subseteq I}{m_F \ne m_G
\mtext{for all} G \subseteq I \mtext{other than} F}.$$
The Scarf complex $\Delta_M$ defines a subcomplex
$\FF_{\Delta_M}$ of the Taylor resolution $\FF_\Delta$.

\proposition{scarfinhull} For any monomial module $M$, the
Scarf complex ${\Delta_M}$ is a subcomplex of the hull
complex $\hull(M)$.

\proof
Let $F = \{ \xx^{\aa_1},\! \ldots,\xx^{\aa_p}\}$
be a face of $\Delta_M$ with $\,m_F = \lcm(F) = \xx^\uu$.
Consider any injective map $\sigma : \{1,\ldots,p\}
\rightarrow \{1,\ldots,n\}$ such that $a_{i,\sigma(i)} = u_i$
for all $i$. Compute the inverse of the $p \times p$-matrix
$(t^{a_{i,\sigma(j)}})$, and let $\vv^\sigma(t)'$ be the sum
of the column vectors of that inverse matrix. By augmenting
the $p$-vector $\vv^\sigma(t)'$ with additional zero coordinates,
we obtain an $n$-vector $\vv^\sigma(t)$ with the following properties:
\item{(i)} $\,  \vv^\sigma(t) \cdot t^{\aa_1} =
 \vv^\sigma(t) \cdot t^{\aa_2} =  \cdots =
\vv^\sigma(t) \cdot t^{\aa_p} \,=\, 1$;
\item{(ii)}
$\, v^\sigma_j(t) \, = \, 0 \,$, for all $j \not\in \image(\sigma)$;
\item{(iii)}
$\, v^\sigma_j(t) \, = \, t^{-u_j} \,+ \,$
{\sl lower order terms in $t$}, for all $j \in \image(\sigma)$.

By taking a convex combination of the vectors $\vv^\sigma(t)$
for all possible injective maps $\sigma$ as above, we obtain a vector
$\vv(t)$ with the following properties:
\item{(iv)} $\,  \vv(t) \cdot t^{\aa_1} =
 \vv(t) \cdot t^{\aa_2} =  \cdots =
\vv(t) \cdot t^{\aa_p} \,=\, 1$;
\item{(v)}
$\, v_j(t) \, = \, c_j \cdot t^{-u_j} \,+ \, $
{\sl lower order terms in $t\,$}
with $c_j > 0 $,  for all $j \in \{1,\ldots,n\}$.

For any $\xx^\bb \in M$ which is not in $F$
there exists an index $\ell$ such that $b_\ell \geq u_\ell + 1$.
This implies $\, \vv(t) \cdot t^\bb
\geq  c_\ell \cdot t^{b_\ell-u_\ell} \,+ \, $ {\sl lower order terms in} $t$,
and therefore $\, \vv(t) \cdot t^\bb > 1\,$ for $t \gg 0$.
We conclude that $F$ defines a face of $P_t$ with inner normal vector
$\vv(t)$. \Box

A binomial first syzygy of $M$ is
called {\it generic} if it has full support, i.e., if no
variable $x_i$ appears with the same exponent in the
corresponding pair of monomial generators. We call
 $M$ {\it generic} if it has a basis of generic binomial
first syzygies. This is a translation-invariant generalization of the
definition of genericity in [BPS].

\lemma{genlem}
 If $M$ is generic, then for any pair of generators $m_i$, $m_j$
either the corresponding binomial first syzygy is generic, or there
exists a third generator $m$ which strictly divides
the least common
multiple of $m_i$ and $m_j$ in all coordinates.

\proof
 Suppose that the syzygy formed by $m_i$ and $m_j$ is not generic, and
induct on the length of a chain of generic syzygies needed to express it.
If the chain has length two, then the middle monomial $m$
divides $\,\lcm(m_i,m_j)$.
Moreover, because the two syzygies involving $m$ are generic, this
division is strict in each variable.
If the chain is longer, then divide it into two
steps. Either each step represents a generic syzygy, and we use the
above argument, or by induction there exists an $m_j$ strictly
dividing the degree of one of these syzygies in all
coordinates, and we are again done.
 \Box

\lemma{notunder}
Let $M$ be a monomial module and $F$ a face of  $\hull(M)$.
For every monomial $m \in M$ there exists a variable
$x_j$ such that $\deg_{x_j}(m) \geq \deg_{x_j}(m_F)$.

\proof
Suppose that $m = \xx^\uu$ strictly divides $m_F$ in each
coordinate. Let $t^{\aa_1},\ldots,t^{\aa_p}$ be the vertices
of $F$ and consider their barycenter
$\,\vv(t) = {1 \over p} \cdot (t^{\aa_1}+ \cdots + t^{\aa_p})\,\in \, F$.
The $j$th coordinate of $\vv(t)$ is a polynomial in $t$
of degree equal to $\,\deg_{x_j}(m_F)$. The $j$th coordinate
of $t^\uu$ is a monomial of strictly lower degree.
Hence $\,t^{\bf u} < \vv(t)\,$ coordinatewise for $t \gg 0$.
Let $\ww$ be a nonzero linear functional which
is nonnegative on  $\RR_+^n $ and whose minimum over
$P_t$ is attained at the face $F$.
Then $\ww \cdot \vv(t)  = \ww \cdot
\aa_1 = \cdots =  \ww \cdot \aa_p $, but
our discussion implies $\ww \cdot t^\uu < \ww \cdot \vv(t) $,
a contradiction. \Box

\theorem{genmin}
 If $M$ is a generic monomial module then $\hull(M)$ coincides with
the Scarf complex $\Delta_M$ of $M$, and the hull resolution
$\FF_{\hull(M)} = \FF_{\Delta_M}$ is minimal.

\proof
Let $F$ be any face of $\hull(M)$ and
$\xx^{\aa_1}, \ldots,\xx^{\aa_p}$ the generators
of $M$ corresponding to the vertices of $F$.
Suppose that $F$ is not a face of $\Delta_M$.  Then either
\item{(i)} $\,
\lcm(\xx^{\aa_1}, \ldots,\xx^{\aa_{i-1}},
\xx^{\aa_{i+1}},\ldots,\xx^{\aa_p}) = m_F \,$ for some $i
\in \{1,\ldots,p\}$, or
\item{(ii)} there exists another generator $\xx^\uu$ of $ M$ which
divides $m_F$ and such that $t^\uu \not\in F$.

Consider first case (i). By \ref{notunder}
applied to $m = {\bf x}^{\aa_i} $ there exists
$x_j$ such that $\deg_{x_j}(\xx^{\aa_i}) = \deg_{x_j}(m_F)$,
and hence $\deg_{x_j}(\xx^{\aa_i}) =
\deg_{x_j}(\xx^{\aa_k})$ for some
$k \not= i$. The first
syzygy between $\xx^{\aa_i}$ and $\xx^{\aa_k}$ is not generic,
and, by \ref{genlem}, there exists a generator $m$ of $M$
which strictly divides $\lcm(\xx^{\aa_i},\xx^{\aa_k})$ in
all coordinates. Since $\lcm(\xx^{\aa_i},\xx^{\aa_k})$
 divides $m_F$, we get a contradiction to \ref{notunder}.

Consider now case (ii).
For any variable $x_j$ there exists $i \in \{1,\ldots,p\}$ such that
$\deg_{x_j}(\xx^{ \aa_i}) = \deg_{x_j}(m_F)
\geq \deg_{x_j}(\xx^\uu)$. If the inequality $\geq$ is an
equality $=$, then the first syzygy between $\xx^\uu$
and $\xx^{\aa_i}$ is not generic, and \ref{genlem} gives
a new monomial generator $m$ which strictly divides $m_F$ in all
coordinates, a contradiction to \ref{notunder}.
Therefore $\geq$ is a strict inequality $>$
for all variables $x_j$. This means that $\xx^\uu$
strictly divides $m_F$ in all coordinates,
again a contradiction to \ref{notunder}.

Hence both (i) and (ii) lead to a contradiction,
and we conclude that every face of $\hull(M)$ is  a face
of $\Delta_M$. This implies $\hull(M) = \Delta_M$ by \ref{scarfinhull}.
The resolution $\FF_{\Delta_M}$ is minimal
because no two faces in $\Delta_M$ have the same degree. \Box

In this paper we are mainly interested in nongeneric
monomial modules for which the hull complex
is typically not simplicial. Nevertheless
the possible combinatorial types of facets
seem to be rather limited. Experimental evidence suggests:

\conjecture{}
Every face of $\hull(M)$ is affinely isomorphic to a subpolytope
of the $(n-1)$-dimensional permutohedron and hence has
at most $\, n \, ! \,$ vertices.

By \ref{exD} it is easy to see that
any subpolytope of the $(n-1)$-dimensional permutohedron
can be realized as the hull complex of suitable monomial ideal.

The following example, found in discussions
with Lev Borisov, shows that the hull complex
of a monomial module need not be locally finite:

\example{lev}
Let $n=3$ and $M$ the monomial module
generated by $x_1^{-1} x_2$ and $\setdef{x_2^i
x_3^{-i}}{i\in\ZZ}$. Then every triangle of the form
$\setthree{x_1^{-1} x_2}{ x_2^i x_3^{-i}}{x_2^{i+1} x_3^{-i-1}}$
is a facet of $\hull(M)$.
In particular, the vertex $x_1^{-1} x_2$ of $\hull(M)$ has
infinite valence.
 \Box

For a generic monomial module $M$ we have the
following important identity
 $$ \, \hull(M_{\preceq \bb}) \quad = \quad \hull(M)_{\preceq \bb}.$$
See equation (5.1) in [BPS]. This identity can fail if $M$ is not generic:

\example{}
 Consider the monomial ideal $M=\idealfour{a^2b}{ac}{b^2\!}{bc^2}$ studied
in \ref{exA} and let $\bb=(2,1,2)$. Then
$\hull(M_{\preceq \bb})$ is a triangle, while $\hull(M)_{\preceq \bb}$
consists of two edges. The vertex $b^2$ of $\hull(M)$
``eclipses'' the facet of $\hull(M_{\preceq \bb})$
 \Box

The hull complex $\hull(M)$ is particularly easy to compute
if $M$ is a squarefree monomial ideal. In this case we have $P_t = P_1$
for all $t$. Moreover, if all square-free generators of $M$ have the
same total degree, then the faces of their convex hull are precisely
the bounded faces of $P_t$.
\ref{thmH} implies the following corrollary.

\corollary{squarefree}
Let $\aa_1,\ldots,\aa_p$ be $0$-$1$-vectors having the same
coordinate sum. Then their boundary complex, consisting of
all faces of the convex polytope $\,P = \conv \{\aa_1,\ldots,\aa_p\}$,
defines a cellular resolution of the ideal
$\,M = \langle  \xx^{\aa_1},\ldots,  \xx^{\aa_p} \rangle$.

\example{exRb} \ref{squarefree} applies
to the Stanley-Reisner ideal of the real projective plane
in \ref{exR}. Here $P$ is  a 5-dimensional polytope with 22
facets, corresponding to the 22 cycles on the $2$-complex
$X$ of length $\le 6$.
Representatives of these three cycle types, and supporting hyperplanes of
the corresponding facets of $P$, are shown on the right in \ref{rpfig}. This
example illustrates how the hull resolution encodes combinatorial information
without making arbitrary choices.
 \Box

\section{} Lattice ideals

\noindent Let $L\subset\ZZ^n$ be a lattice. In this section
we study (cellular) resolutions of the lattice module $M_L$ and of the
lattice ideal $I_L$. Let $S[L]$ be the group algebra of $L$
over $S$. We realize $S[L]$ as the subalgebra of
$\,k[x_1,\ldots,x_n, z_1^{\pm 1}\!\!,\, \dots ,\, z_n^{\pm 1}]
\,$ spanned by all monomials
$\xx^\aa \zz^\bb$ where $\aa \in \NN^n$ and $\bb \in L$. Note that
$\,S \,= \,S[L]/ \langle \,\zz^\aa - 1 \,\,| \,\, \aa \in L\, \rangle $.

\lemma{lem1}
The lattice module $M_L$ is an $S[L]$-module,
and $\,M_L \otimes_{S[L]} S \,= \, S/I_L$.

\proof
The $k$-linear map
$\,\phi : S[L] \rightarrow M_L , \, \xx^\aa \zz^\bb
\mapsto \xx^{\aa + \bb}   \,$
defines the structure of an $S[L]$-module on $M_L$.
Its kernel $\,\ker(\phi) \,$ is the ideal in $S[L]$ generated by all
binomials $\,\xx^\uu - \xx^\vv \zz^{\uu-\vv}$
where $\uu,\vv \in \NN^n$ and $\uu - \vv \in L$.
Clearly, we obtain $I_L$ from $\ker(\phi)$
by setting all $\zz$-variables to $1$, and hence
$\,(S[L]/\ker(\phi) )\otimes_{S[L]} S \,= \, S/I_L$. \Box

We define a $\ZZ^n$-grading on $S[L]$ via
$\,\deg(\xx^\aa \zz^\bb) = \aa + \bb$.
Let ${\cal A}$ be the category of $\ZZ^n$-graded $S[L]$-modules,
where the  morphisms  are $\ZZ^n$-graded $S[L]$-module homomorphisms
of degree $\zero$. The polynomial ring
$S = k[x_1,\ldots,x_n]$ is graded by the quotient group
$\ZZ^n/L$ via $\,\deg(\xx^\aa) = \aa + L$. Let ${\cal B}$
be the category of $\ZZ^n/L$-graded $S$-modules, where the
morphisms  are  $\ZZ^n/L$-graded $S$-module homomorphisms
of degree $\zero$.     Clearly, $M_L$ is an object in ${\cal A}$,
and $\,M_L \otimes_{S[L]} S \,= \, S/I_L \,$ is an object in ${\cal B}$.

\theorem{thmU}
The categories  ${\cal A}$ and ${\cal B}$ are equivalent.

\proof
 Define a functor $\pi:{\cal A} \rightarrow {\cal B}$ by the rule
$\,\pi(M) :=  M \otimes_{S[L]} S$. This functor  weakens
the $\ZZ^n$-grading of objects in ${\cal A}$ to a $\ZZ^n/L$-grading.
The  properties of $\pi$ cannot be deduced from the
tensor product alone, which is poorly behaved when applied to
arbitrary $S[L]$-modules; e.g., $S$ is not a flat $S[L]$-module.
Further, the categories
${\cal A}$ and ${\cal B}$ are {\sl not isomorphic};
we are only claiming that  they are {\sl equivalent}.

We apply condition iii) of [Mac, \S IV.4, Theorem 1]: It is enough to prove
that $\pi$ is full and faithful, and that each object $N\in{\cal B}$ is
isomorphic to $\pi(M)$ for some object $M\in{\cal A}$. To prove that $\pi$
is full and faithful, we show that for any two modules $M$, $M'\in {\cal
A}$ it induces an identification $\Hom_{\cal
A}(M,M') = \Hom_{\cal B}(\pi(M),\pi(M'))$.

Because each module $M\in{\cal A}$ is $\ZZ^n$-graded, the lattice
$L\subset S[L]$ acts on $M$ as a group of automorphisms, i.e. the
multiplication maps $\zz^\bb: M_\aa \rightarrow M_{\aa+\bb}$ are
isomorphisms of $k$-vector spaces for each $\bb \in L$, compatible with
multiplication by each $x_i$. For each $\alpha\in\ZZ^n/L$, the functor
$\pi$ identifies the spaces $M_\aa$ for $\aa\in\alpha$ as the single space
$\pi(M)_\alpha$.
A morphism $f:M\rightarrow M'$ in ${\cal A}$ is a collection of $k$-linear
maps $f_\aa: M_\aa\rightarrow M'_\aa$, compatible with the action by $L$
and with multiplication by each $x_i$. A morphism $g:\pi(M)\rightarrow
\pi(M')$ in ${\cal B}$ is a collection of $k$-linear maps
$g_\alpha:\pi(M)_\alpha\rightarrow \pi(M')_\alpha$, compatible with
multiplication by each $x_i$. For each $\alpha\in\ZZ^n/L$, the functor
$\pi$ identifies the maps $f_\aa$ for $\aa\in\alpha$ as the single map
$\pi(f)_\alpha$.

It is clear from this discussion that $\pi$ takes
distinct morphisms to distinct morphisms. Given a morphism $g\in\Hom_{\cal
B}(\pi(M),\pi(M'))$, define
a morphism $f\in\Hom_{\cal A}(M,M')$ by the rule $f_\aa=g_\alpha$ for
$\aa\in\alpha$. We have $\pi(f)=g$, establishing the desired identification
of Hom-sets. Hence $\pi$ is full and faithful.

Finally, let $N = \dsum_{ \alpha \in \ZZ^n/L} N_\alpha$ be any object in
${\cal B}$. We define an object  $M = \oplus_{\aa \in \ZZ^n} M_\aa$ in
${\cal A}$ by setting $M_\aa := N_\alpha$ for each $\aa \in \alpha$, by
lifting each multiplication map $x_i:N_\alpha \rightarrow N_{\alpha+\ee_i}$
to maps $x_i:M_\aa \rightarrow M_{\aa+\ee_i}$ for $\aa \in \alpha$, and
by letting $\zz^\bb$ act on $M$ as the identity map from $M_\aa $ to
$M_{\aa+\bb}$ for $\bb \in L$. The module $M$ satisfies $\pi(M) = N$,
showing that $\pi$ is an equivalence of categories.
 \Box

\ref{thmU} allows us to resolve the lattice module $M_L\in
{\cal A}$ in order
to resolve the quotient ring $\pi(M_L) = S/I_L\in {\cal B}$,
and conversely.

\corollary{corU}
A $\ZZ^n$-graded complex  of free $S[L]$-modules,
 $$ C : \qquad \;\cdots \;
\rightarrowbox{8pt}{$f_2$}\; S[L]^{\beta_1}
\rightarrowbox{8pt}{$f_1$}\; S[L]^{\beta_0}
\;\rightarrowbox{8pt}{$f_0$}\; S[L] \;\rightarrow \; M_L \;\rightarrow \; 0, $$
 is a (minimal) free resolution of $M_L$ if and only if its image
 $$  \pi(C)\,: \, \;\cdots \;
\rightarrowbox{8pt}{$\pi(f_2)$}\; S^{\beta_1}
\rightarrowbox{8pt}{$\pi(f_1)$}\; S^{\beta_0}
\;\rightarrowbox{8pt}{$\pi(f_0)$}\; S
\;\rightarrow \; S/I_L
\;\rightarrow \; 0 ,$$
is a (minimal) $\ZZ^n/L$-graded  resolution of $S/I_L$
by free $S$-modules.

\proof
 This follows immediately from \ref{lem1} and \ref{thmU}.
 \Box

Since $S[L]$ is a free $S$-module, every resolution $C$
as in the previous corollary gives rise to
a resolution of $M_L$ as a $\ZZ^n$-graded $S$-module.  We demonstrate
in an example how resolutions of $M_L$ over $S$
are derived from resolutions of $S/I_L$ over $S$.

\example{exU}
 Let $S=k[x_1,x_2,x_3]$ and $L=\ker \bmatrix{2pt}{1 & 1 &  1 \cr}
 \subset \ZZ^3$.
Then $\ZZ^3 /L \simeq \ZZ$,   $I_L = \ideal{x_1-x_2,x_2-x_3}$,
and $M_L$ is the module generated by all
monomials of the form $\ x_1^i x_2^j x_3^{-i-j}$.
The ring $S/I_L$ is resolved by the Koszul complex
 $$ 0 \longrightarrow
 S(-2) \rightarrowmat{4pt}{4pt}{\!\!x_2 - x_3 \cr x_2 - x_1 \cr}
 S(-1)^2 \rightarrowmat{4pt}{6pt}{x_1 - x_2 & x_2 - x_3 \cr}
 S \longrightarrow S/I_L . $$
This is a $\ZZ^3/L$-graded complex of free $S$-modules.
An inverse image under $\pi$ equals
 $$ \eqalign{
 0  \longrightarrow
 S[L]\bigl(-(1,1,0)\bigr)
& \rightarrowmat{4pt}{4pt}{\!\!x_2 - x_3 z_2  z_3^{-1}\cr
 x_2- x_1 z_2 z_1^{-1}  \cr}
 S[L]\bigl(-(1,0,0)\bigr)  \oplus
 S[L]\bigl(-(0,1,0)\bigr) \cr
 & \qquad \qquad
\rightarrowmat{4pt}{6pt}{x_1 - x_2 z_1 z_2^{-1} & x_2 - x_3 z_2
z_3^{-1} \cr}  S[L] \longrightarrow M_L .\cr}  $$
Writing each term as a direct sum of free $S$-modules,
for instance, $\, S[L]\bigl(-(1,1,0)\bigr) =
\oplus_{i+j+k=2} S \bigl(-(i,j,k) \bigr)$, we
get a $\ZZ^3$-graded minimal free resolution of $M_L$ over $S$:
$$
0 \,\, \rightarrow
\dsum_{i+j+k=2}\! S \bigl(-(i,j,k) \bigr)
 \,\, \rightarrow
\dsum_{i+j+k=1}\! S \bigl(-(i,j,k) \bigr)^2
 \,\, \rightarrow
\dsum_{i+j+k=0}\! S \bigl(-(i,j,k) \bigr)
 \rightarrow
M_L. \Box
$$

\vskip .2cm

Our goal is to define and study
cellular resolutions of the lattice ideal $I_L$.
Let $X$ be a $\ZZ^n$-graded cell complex whose vertices
are the generators of $M_L$.
Each cell $F \in X$ is identified with its set of vertices,
regarded as  a subset of $L$.
The cell complex $X$ is called {\it equivariant}
\ if $\,F \in X$ and $\bb \in L$ implies that $F + \bb \,\in \, X$,
and if the incidence function satisfies
$\,\varepsilon(F,F') = \varepsilon(F+\bb,F'+\bb)$
for all $\bb \in L$.

\lemma{ }
If $X$ is an equivariant $\ZZ^n$-graded cell complex on $M_L$
then the cellular complex $\FF_X$ has the structure of a
$\ZZ^n$-graded complex of free $S[L]$-modules.

\proof
The group $L$ acts on the faces of $X$.
Let $X/L$ denote the set of orbits. For each orbit ${\cal F}
\in X/L $ we select a distinguished representative $F \in {\cal F}$,
and we write $\Rep(X/L)$ for the set of representatives.
The following map is an
isomorphism of $\ZZ^n$-graded $S$-modules, which
defines the structure of a free $S[L]$-module on $\FF_X$:
$$  \dsum_{F \in \Rep(X/L)} \!\!\!\! S[L] \cdot e_F \quad \simeq \quad
  \dsum_{F \in X }S \cdot e_F \, \,
\,\,\, = \,\,\, \FF_X , \quad \,
\zz^\bb \cdot e_F \,\, \mapsto \,\,  e_{F + \bb}\ .$$
The differential $\partial $ on $\FF_X$
is compatible with the $S[L]$-action
on $\FF_X$ because the incidence
function is $L$-invariant. For each $F \in \Rep(X/L)$
and $\bb \in L$ we have
$$ \eqalign{
 \partial(\zz^\bb \cdot e_F ) \quad  &= \quad
 \partial ( e_{F+\bb}) \quad
= \quad \sum_{F'\in X,\, F'\ne\emptyset} \;
  \varepsilon(F \! + \! \bb,F' \! + \!\bb) \; {m_{F \!+\!  \bb}
\over m_{F'\! +\! \bb}}  \; e_{F'+\bb}  \cr
& = \,\,  \sum_{F'\in X,\, F'\ne\emptyset} \; \! \!
\varepsilon(F,F') \; {m_{F} \over m_{F'}} \; \zz^\bb  \cdot e_{F'}
\quad = \quad
\zz^\bb \cdot \partial(  e_F ) .  \cr} $$
Clearly, the differential $\partial$ is homogeneous of degree $0$,
which proves the claim. \Box

\corollary{}
If $X$ is an equivariant $\ZZ^n$-graded cell complex on $M_L$
then the cellular complex $\FF_X$ is exact over $S$ if and only if
it is exact over $S[L]$.

\proof The $\ZZ^n$-graded components of $\FF_X$ are
complexes of $k$-vector spaces which are independent of our
interpretation of $\FF_X$ as an $S$-module or  $S[L]$-module. \Box

If $X$ is an equivariant $\ZZ^n$-graded cell complex on $M_L$
such that $\FF_X$ is exact, then
we call $\FF_X$ an {\it equivariant cellular resolution} of $M_L$.

\corollary{cor17}
If $\FF_X$ is an equivariant cellular (minimal) resolution of $M_L$
then $\, \pi(\FF_X)\,$ is a (minimal) resolution of $S/I_L$
by $\ZZ^n/L$-graded free $S$-modules.
 \Box

We call $\,\pi(\FF_X)\,$ a
{\it cellular resolution} of the lattice ideal $I_L$.
Let $Q$ be an order ideal in the quotient poset $\NN^n/L$.
Then $Q + L$ is an order ideal in $\NN^n +L$, and
the restriction $\,\FF_{X_{Q+L}}\,$ is a complex
of $\ZZ^n$-graded  free $S[L]$-modules. We set
$\,\pi(\FF_X)_Q \,:= \,\pi(\FF_{X_{Q+L}})$. This
is a complex of $\ZZ^n/L$-graded free $S$-modules.
\ref{corK} implies

\proposition{restr}
If $\pi(\FF_X)$ is a cellular resolution of $I_L$ and
$Q$ is an order ideal in $\NN^n/L$ which
contains all Betti degrees then
$\,\pi(\FF_{X})_Q\,$ is a cellular resolution of $I_L$.

In what follows we shall study two particular
cellular resolutions of $I_L$.

\theorem{nicethm} The Taylor complex $\Delta$ on $M_L$
and the hull complex $\hull(M_L)$ are equivariant.
They define cellular resolutions
$\,\pi(\FF_\Delta)\,$  and $\,\pi(\FF_{\hull(M_L)}) \,$ of $I_L$.

\proof
The Taylor complex $\Delta$
consists of all finite subsets of generators of $M_L$.
It has an obvious $L$-action. The hull complex also has an $L$-action: if
$\,F = \conv \bigl( \set{ t^{\aa_1},\ldots,t^{\aa_s} } \bigr) \,$
is a face of $\hull(M_L)$ then
 $\,\zz^b \cdot F =
\conv \bigl( \set{ t^{\aa_1+\bb},\ldots,t^{\aa_s+\bb} } \bigr) \,$
is also a face of $\hull(M_L)$ for all $\bb \in L$.
In both cases the incidence function $\varepsilon$ is defined uniquely
by the ordering of the elements in $L$. To ensure that
$\varepsilon$ is $L$-invariant, we
fix an ordering which is $L$-invariant; for instance,
order the elements of $L$ by the value of an
$\RR$-linear functional whose coordinates are $\QQ$-linearly independent.

Both $\pi(\FF_\Delta)$ and $\pi(\FF_{\hull(M_L)})$
are cellular resolutions of $I_L$
by \ref{cor17}.
\Box

The {\it Taylor resolution} $\pi(\FF_\Delta)$ of $I_L$ has
the following explicit description.
For $\alpha \in \NN^n/L$ let $\fiber(\alpha)$ denote the
(finite) set of all monomials $\xx^\bb$ with $\bb \in \alpha$.
Thus $S_\alpha = k \cdot \fiber(\alpha)$.
Let $E_i(\alpha)$  be the collection of
all $i$-element subsets $I$ of $\fiber(\alpha)$
whose greatest common divisor $\gcd(I)$ equals $1$.
For $I \in E_i(\alpha)$ set $\deg(I) := \alpha$.

\proposition{explicit}
The Taylor resolution $\pi(\FF_\Delta)$ of
a lattice ideal $I_L$ is isomorphic to the
$\ZZ^n/L$-graded free $S$-module
$\,\, \dsum_{\alpha \in \NN^n/L} S \cdot E_i(\alpha)\,$
with the differential
$$ \partial(I)\quad = \quad \sum_{m \in I} \sign(m,I) \cdot
\gcd(I \backslash \set{m})\cdot [I \backslash \set{m}]. \eqno (3.1) $$
In this formula, $\,[I  \backslash \set{m}] $ denotes the
element of
$ \,E_{i-1}\bigl( \alpha - \deg(\gcd(I \backslash \set{m}))\bigr)$
which is obtained from $I \backslash \set{m}$
 by removing the common factor $\,\gcd(I  \backslash \set{m})$.

\proof
For $\bb \in \ZZ^n$ let $F_i(\bb)$ denote the
collection of $i$-element subsets of generators of $M_L$
whose least common multiple equals $\bb$.
For $J \in F_i(\bb)$ we have $\lcm(J) = \xx^\bb$. The
Taylor resolution $\FF_\Delta$ of $M_L$ equals
$\, \dsum_{\bb \in \NN^n + L} S \cdot F_i(\bb)\,$
with differential
$$ \partial(J) \quad = \quad \sum_{m \in J} \sign(m,J) \cdot
{\lcm(J) \over \lcm( J \backslash \set{m})} \cdot J \backslash \set{m}.
\eqno (3.2) $$
There is a natural  bijection between $F_i(\bb)$ and $E_i(\bb+L)$,
namely, $\,J \, \mapsto  \, \set{ \xx^\bb / \xx^\cc
\,\,|\,\,\xx^\cc \in J} \, = \, I $.
Under this bijection we have
$\,{\xx^\bb \over \lcm( J \backslash \set{m})}
= \gcd(I \backslash \set{m})$.
The functor $\pi$ identifies each $F_i(\bb)$ with
$E_i(\bb+L)$ and it takes  (3.2) to (3.1). \Box

\corollary{}
Let $Q$ be an order ideal in $\NN^n/L$ which
contains all Betti degrees.
Then  $\,\pi(\FF_\Delta)_Q
= \dsum_{\alpha \in Q} S E_i(\alpha)\,$
with differential (3.1) is a cellular resolution of $I_L$.

\proof
This follows from \ref{restr}, \ref{nicethm}
and \ref{explicit}.
\Box

\example{} {\sl (Generic lattice ideals) }
The lattice module $M_L$ is generic (in the sense of \S 2)
if and only if the  ideal $I_L$ is generated by binomials with full support.
Suppose that this holds. It was shown in [PS] that
the Betti degrees of $I_L$ form an order ideal $Q$
in $\NN^n/L$. \ref{genmin} and \ref{restr} imply that
the resolution $\, \pi(\FF_\Delta)_Q $ is minimal and coincides with
the hull resolution $\pi(\FF_{\hull(M_L)})$.
\Box

\vskip .2cm

The remainder of this section is devoted to the hull resolution
of $I_L$. We next show that the hull complex
$\hull(M_L)$ is locally finite. This fact is nontrivial,
in view of \ref{lev}. It will imply that the hull
resolution has finite rank over $S$.

Write each vector $\aa \in L \subset \ZZ^n$ as difference
$\aa = \aa^+ - \aa^-$ of two nonnegative vectors
with disjoint support. A nonzero vector $\aa \in L$
is called {\it primitive} if there is no vector
$\bb \in L \backslash \set{ \aa, \zero }$ such that
$\bb^+ \leq \aa^+$ and $\bb^- \leq \aa^-$.
The set of primitive vectors is known to be finite
[St, Theorem 4.7]. The set of binomials
 $\,\xx^{\aa^+} - \xx^{\aa^-}\,$ were $\aa$ runs over
all primitive vectors in $L$ is called
the {\it Graver basis} of the ideal $I_L$.
The  Graver basis contains the universal
Gr\"obner basis of $I_L$ [St, Lemma 4.6].

\lemma{graver}
If $\set{\zero,\aa}$ is an edge of $\hull(M_L)$ then
$\aa$ is a primitive vector in $L$.

\proof
Suppose that $\aa= (a_1,\ldots,a_n) $
is a vector in $L$ which is not primitive,
and choose $\,\bb  = (b_1,\ldots,b_n)
\in L \backslash \set{ \aa, \zero}$ such that
$\bb^+ \leq \aa^+$ and $\bb^- \leq \aa^-$.
This implies $\,t^{b_i} + t^{a_i-b_i} \le 1 + t^{a_i} \,$ for
$t \gg 0$ and $i \in \set{1,\ldots,n}$.
In other words, the vector $\,t^{\bb} + t^{\aa-\bb}\,$
is componentwise smaller or equal to the vector
$\, t^\zero + t^\aa $. We conclude that the
midpoint of the segment $ \, \conv \set{ t^\zero , t^\aa}\, $
lies in $\,\conv \set{t^{\bb} , t^{\aa-\bb}}  + \RR_+^n$,
and hence $\,\conv \set{ t^\zero , t^\aa}$ is not an edge of
the polyhedron
$\,P_t \, = \, \conv \set{\, t^\cc \,: \, \cc \in L} \,+ \, \RR_+^n$.
 \Box

\theorem{ }
The hull resolution $\pi(\FF_{\hull(M_L)})$ is finite as an $S$-module.

\proof
By \ref{graver} the vertex $\zero$ of $\hull(M_L)$ lies in only
finitely many edges. It follows that $\zero$ lies in only finitely
many faces of $\hull(M_L)$. The lattice $L$ acts transitively
on the vertices of $\hull(M_L)$, and hence every face of $\hull(M_L)$
is $L$-equivalent to a face containing $\zero$.
The faces containing $\zero$ generate
$\FF_{\hull(M_L)}$ as an $S[L]$-module, and hence they
generate $\pi(\FF_{\hull(M_L)})$ as an $S$-module. \Box

A minimal free resolution of a lattice ideal $I_L$
generally does not respect symmetries, but
the hull resolution does.
The following example illustrates this point.

 \example{exE}
 {\sl (The hypersimplicial complex as  a hull resolution)}
 \hfill \break The lattice $\,L \,=\,  \ker_\ZZ \pmatrix{ 1 \! & \!1 &
\cdots & 1 }\,$ in
$\ZZ^n$  defines the toric ideal
$$ I_L \quad = \quad
\langle \, x_i - x_j \,\,: \,\, 1 \leq i < j \leq n \, \rangle .$$
The minimal free resolution of $I_L$ is
the Koszul complex on $n-1$ of the generators $x_i-x_j$.
Such a minimal resolution does not respect the
action of the symmetric group $S_n$ on $I_L$.
The hull resolution is the
Eagon-Northcott complex of the matrix
{\smallmath $\,\bmatrix{2pt}{1\;\, & 1\;\, & \cdots & 1\;\, \cr
        x_1 & x_2 & \cdots & x_n \cr}$}.
This resolution  is not minimal but it retains the
$S_n$-symmetry of $I_L$.
It coincides with the {\sl hypersimplicial complex} studied
by Gel'fand and MacPherson in [GM, \S 2.1.3].
The basis vectors
of the hypersimplicial complex are denoted $\Delta_{\ell}^I$
where $I$ is a subset of $\set{1,2,\ldots,n}$ with $|I| \geq 2$
and $\ell$ is an integer with $ 1 \leq \ell \leq |I| - 1$.
We have $\,\Delta_1^{\set{i,j}} \mapsto  x_i - x_j\,$ and the
higher differentials act as
$$ \Delta_\ell^I \quad \mapsto \quad
\sum_{i \in I} \sign(i,I) \cdot x_i \cdot \Delta^{I \setminus \!
\set{i}}_{\ell-1} \, - \, \sum_{i \in I} \sign(i,I) \cdot \Delta^{I \setminus
\! \set{i}}_{\ell}, $$ where  the first sum is zero if $\ell=1$
and the second sum is zero if $\ell = \abs{I}-1$.
\Box

\remark{curious}
Our study suggests a {\sl curious duality} of toric
varieties, under which the coordinate ring of the primal variety is
resolved by a discrete subgroup of the dual variety.
More precisely, the hull resolution of $I_L$  is gotten by taking the
convex hull in $\RR^n$ of the points $t^\aa$ for $\aa \in L$.
The Zariski closure of these points (as $t$ varies)
is itself an affine toric variety, namely, it is the
variety defined by the lattice ideal
$\,I_{L^\perp}$ where
$L^\perp$ is the lattice dual to $L$
under the standard inner product on $\ZZ^n$.

For instance, in \ref{exE} the primal toric
variety is the line $\,(t , t , \ldots , t)$
and  the dual toric variety is the
hypersurface $x_1 x_2 \cdots x_n = 1$.
That hypersurface forms a group under
coordinatewise multiplication,
and we are taking the convex hull of a
discrete subgroup to resolve the coordinate ring
of the line $\,(t , t , \ldots , t)$. \Box

\example{exV} {\sl (The rational normal quartic curve in $P^4$)}

\vskip .1cm
\noindent
Let $L = \ker_\ZZ {\smallmath\bmatrix{2pt}{
0 & 1 & 2 & 3 & 4 \cr
4 & 3 & 2 & 1 & 0 \cr}}$.
The minimal free resolution of the lattice ideal $I_L$ looks like
$ 0 \rightarrow S^3  \rightarrow S^8 \rightarrow S^6  \rightarrow I_L $.
The primal toric variety in the sense of  \ref{curious}
is a curve in $P^4$ and the dual toric variety
is the embedding of the $3$-torus into affine $5$-space
given by the equations $
\, x_2  x_3^2 x_4^3 x_5^4 =   x_1^4 x_2^3 x_3^2 x_4 = 1$.
Here the hull complex $\hull(M_L)$ is simplicial,
and the hull resolution  of $I_L$ has the format
$ 0  \rightarrow S^4  \rightarrow S^{16}
 \rightarrow S^{20}  \rightarrow S^9
 \rightarrow I_L $.
The nine classes of edges in $\hull(M_L)$
are the seven quadratic binomials in $I_L$
and the two cubic  binomials
$\,x_3 x_4^2 - x_1 x_5^2 , \, x_2^2 x_3 - x_1^2 x_5 $.

\vskip 1.2cm

\noindent
{\bf Acknowledgements. }
We thank Lev Borisov, David Eisenbud, Irena Peeva, Sorin
Popescu, and Herb Scarf for helpful conversations. Dave Bayer and Bernd
Sturmfels are partially supported by the National Science Foundation.
Bernd Sturmfels is also supported by the David and Lucille Packard
Foundation and a 1997/98 visiting position at
the Research Institute for Mathematical Sciences
of Kyoto University.

\bigskip \bigskip

\references

\itemitem{[AB]} R.~Adin and D.~Blanc,
Resolutions of associative and Lie algebras,
preprint, 1997.

\itemitem{[BHS]} I.~Barany, R.~Howe, H.~Scarf:
The complex of maximal lattice free simplices,
{\sl Mathematical Programming} {\bf 66} (1994) Ser.~A, 273--281.

\itemitem{[BPS]} D.~Bayer, I.~Peeva and B.~Sturmfels,
Monomial resolutions, preprint, 1996.

\itemitem{[BLSWZ]} A.~Bj\"orner, M.~Las~Vergnas,
B.~Sturmfels, N.~White and G.~Ziegler,
{\sl Oriented Matroids}, Cambridge University Press, 1993.

\itemitem{[BH]} W.~Bruns and J.~Herzog, {\sl Cohen-Macaulay Rings},
Cambridge University Press, 1993.

\itemitem{[BH1]} W.~Bruns and J.~Herzog,
On multigraded resolutions, {\sl Math.~Proc.~Cambridge Philos.~Soc.}
{\bf 118} (1995) 245--257.

\itemitem{[GM]} I.~M.~Gel'fand and R.~D.~MacPherson:
Geometry in Grassmannians and a generalization of
the dilogarithm, {\sl Advances in Math.}
{\bf 44} (1982), 279--312.

\itemitem{[Ho]} M. Hochster,
Cohen-Macaulay rings, combinatorics and simplicial complexes,
in {\sl Ring Theory II}, eds.~B.R.~McDonald and R.~Morris,
Lecture Notes in Pure and Appl.~Math.~{\bf 26},
Dekker, New York, (1977), 171--223.

\itemitem{[KS]} M.~Kapranov and M.~Saito,
Hidden Stasheff polytopes in algebraic K-theory and the space of
Morse functions, preprint, 1997,
paper \# 192 in {\tt http://www.math.uiuc.edu/K-theory/}.

\itemitem{[Mac]} S. MacLane,
{\sl Categories for the Working Mathematician},
Graduate Texts in Mathematics, No.~5,
Springer-Verlag, New York, 1971.

\itemitem{[PS]} I.~Peeva and B.~Sturmfels,
Generic lattice ideals, to appear in {\sl Journal of the American
Math.~Soc.}

\itemitem{[Ros]} I.~Z.~Rosenknop, Polynomial ideals that are generated by
monomials (Russian), {\sl Moskov. Oblast. Ped. Inst. Uw cen Zap.} {\bf 282}
(1970), 151-159.

\itemitem{[Stu]} B.~Sturmfels,
{\sl Gr\"obner Bases and Convex Polytopes},
AMS University Lecture Series,  Vol. 8, Providence RI, 1995.

\itemitem{[Stu2]} B.~Sturmfels,
The co-Scarf resolution, to appear in
{\sl Commutative Algebra and Algebraic Geometry},
Proceedings Hanoi 1996, eds.~D.~Eisenbud and N.V.~Trung,
Springer Verlag.

\itemitem{[Tay]} D.~Taylor,
{\sl Ideals Generated by Monomials in an $R$-Sequence},
Ph.~D.~thesis, University of Chicago, 1966.

\itemitem{[Zie]} G.~Ziegler, {\sl Lectures on Polytopes},
Springer, New York, 1995.

\vskip 1.8cm

\noindent
 Dave Bayer, Department of Mathematics,
Barnard College, Columbia University,
New York, NY 10027, USA,
{\tt bayer@math.columbia.edu}.

\vskip .4cm

\noindent
Bernd Sturmfels,
Department of Mathematics,
University of California,
Berkeley, CA 94720, USA,
{\tt bernd@math.berkeley.edu}.
\bye